\documentclass[12pt]{article}

\usepackage[a4paper,margin=2.4cm]{geometry}
\usepackage{amsmath,amssymb,amsthm}
\usepackage{booktabs}
\usepackage[hidelinks]{hyperref}
\theoremstyle{plain}
\newtheorem{theorem}{Theorem}[section]
\newtheorem{proposition}[theorem]{Proposition}

\newtheorem{corollary}[theorem]{Corollary}
\theoremstyle{definition}
\newtheorem{definition}[theorem]{Definition}
\newtheorem{example}[theorem]{Example}
\newtheorem{hypothesis}[theorem]{Hypothesis}
\theoremstyle{remark}
\newtheorem{remark}[theorem]{Remark}

\renewcommand{\Re}{\operatorname{Re}}
\renewcommand{\Im}{\operatorname{Im}}

\title{Zeta Functions of Polygonal Numbers:\\
Analytic Continuation, Poles and Zero Scattering}
\author{Jan Luegering}
\date{November 16th 2026}

\begin{document}
\maketitle

\begin{abstract}
We study the Zeta functions associated with polygonal numbers, focusing on their meromorphic continuation, poles, residues, and trivial zeros. We then investigate the distribution of their nontrivial zeros. Under explicit hypotheses, we derive an asymptotic formula for the mean real part of zeros in long height windows. The formula predicts convergence of this mean to $1/4$ as the height tends to infinity. We also formulate a stronger concentration conjecture and examine numerical evidence for it. 
\end{abstract}

\newpage

\tableofcontents

\newpage

\section{Definition}\label{sec:definition}
\begin{definition}[Polygonal numbers]
For an integer $k>2$ the $n$-th
$k$-gonal number is
\begin{equation}\label{eq:Pk}
  P_k(n)\;=\;\frac{(k-2)\,n^{2}-(k-4)\,n}{2}.
\end{equation}
\end{definition}

The cases $k=3,4,5$ give the triangular numbers $1,3,6,10,\dots$, the
squares $1,4,9,16,\dots$, and the pentagonal numbers $1,5,12,22,\dots$,
respectively.

It is convenient to introduce two constants,
\begin{equation}\label{eq:ab}
  a\;=\;\frac{k-2}{2},
  \qquad
  b\;=\;\frac{k-4}{k-2}\;=\;1-\frac{2}{k-2}.
\end{equation}
Then $P_k(n)$ can be written in three equivalent ways:
\begin{equation}\label{eq:Pk-forms}
  P_k(n)\;=\;a\,n\,(n-b)
  \;=\;n+(k-2)\binom{n}{2}
  \;=\;a\Bigl[\bigl(n-\tfrac b2\bigr)^{2}-\tfrac{b^{2}}{4}\Bigr].
\end{equation}
We will use the identities
\begin{equation}\label{eq:first-term}
  P_k(1)=1\quad\text{for every }k,
  \qquad
  a\,(1-b)=1 .
\end{equation}

\begin{definition}[Zeta function of the polygonal numbers]
\begin{equation}\label{eq:Zk}
  Z_k(s)\;=\;\sum_{n=1}^{\infty}P_k(n)^{-s},\qquad s=\sigma+it .
\end{equation}
\end{definition}

Although polygonal numbers correspond to integer values of $k$, it is
useful in the analytic arguments to allow real $k>2$ as well. Since
$P_k(n)\sim a n^2$, the defining series converges absolutely in the
half-plane $\sigma > \tfrac12$.

\begin{example}[Two classical cases]
For the excluded case $k=2$ we have $P_2(n)=n$, so $Z_2(s)=\zeta(s)$.
For $k=4$ we have $P_4(n)=n^{2}$, so $Z_4(s)=\zeta(2s)$.
In these two cases the functions inherit the familiar Euler product and
functional equation from the Riemann zeta function.
\end{example}

\paragraph{Factorisation.}
Writing $P_k(n)=a\,n(n-b)$ gives
\begin{equation}\label{eq:factor}
  Z_k(s)\;=\;a^{-s}\,H_b(s),\qquad
  H_b(s)\;=\;\sum_{n\ge1}\bigl(n(n-b)\bigr)^{-s}.
\end{equation}
The factor $a^{-s}$ never vanishes, so $Z_k$ and $H_b$ have exactly the same
zeros.
It does, however, change the \emph{size} of the function: $|a^{-s}|=a^{-\sigma}$.
The factor $a^{-s}$ will contribute the term $\ln a$ to the asymptotic
formula derived in Section~
ef{sec:law}.

\newpage

\section{Analytic continuation}\label{sec:continuation}
Let $\zeta(w,x)=\sum_{m\ge0}(m+x)^{-w}$ be the Hurwitz zeta function, and let
$(s)_j=s(s+1)\cdots(s+j-1)$, with $(s)_0=1$.

\begin{theorem}[Continuation, all integers $k>2$]\label{thm:continuation}
Put $x=1-\tfrac b2$. Then
\begin{equation}\label{eq:continuation}
  Z_k(s)\;=\;a^{-s}\sum_{j=0}^{\infty}\frac{(s)_j}{j!}
  \Bigl(\frac b2\Bigr)^{2j}\,\zeta(2s+2j,\;x).
\end{equation}
The series converges locally uniformly on compact subsets away from the
poles of its summands. It therefore defines a meromorphic function on
$\mathbb C$.
\end{theorem}

\begin{proof}
Put $m=n-1$ and use the third form in \eqref{eq:Pk-forms}:
\[
  P_k(n)^{-s}
  =a^{-s}(m+x)^{-2s}\Bigl(1-\frac{(b/2)^{2}}{(m+x)^{2}}\Bigr)^{-s}.
\]
Expand $(1-y)^{-s}=\sum_{j\ge0}\frac{(s)_j}{j!}\,y^{j}$ with
$y=\frac{(b/2)^{2}}{(m+x)^{2}}$.
Since $m+x\ge x=1-\tfrac b2>\tfrac{|b|}{2}$ for every $k>2$ (recall $b<1$),
we have $y<1$ and the expansion converges.
Summing over $m\ge0$ turns $\sum_m(m+x)^{-2s-2j}$ into
$\zeta(2s+2j,x)$.
For $\sigma>\tfrac12$ all sums converge absolutely, so they may be
interchanged.
On any compact set of $s$-values, $(s)_j/j!$ grows at most polynomially
in $j$, uniformly in $s$, while
$\zeta(2s+2j,x)=x^{-2s-2j}(1+O(q^j))$ for some $q<1$ (uniformly on that
compact set, for sufficiently large $j$). Thus the tail is bounded by a
polynomial in $j$ times the geometric factor
$\bigl(\tfrac{|b|}{2x}\bigr)^{2j}$. The series therefore converges
locally uniformly away from the displayed term poles, which proves the
claimed meromorphic continuation.
\end{proof}

\begin{example}[The three simplest cases]
\begin{itemize}
\item $k=4$: $b=0$, $a=1$, only $j=0$ survives and $Z_4(s)=\zeta(2s)$.
\item $k=3$: $b=-1$, $a=\tfrac12$, $x=\tfrac32$, so
\begin{equation}\label{eq:Z3}
  Z_3(s)\;=\;2^{s}\sum_{j\ge0}\frac{(s)_j}{j!}\,4^{-j}\,\zeta\bigl(2s+2j,\tfrac32\bigr).
\end{equation}
\item $k=5$: $b=\tfrac13$, $a=\tfrac32$, $x=\tfrac56$.
\end{itemize}
\end{example}

\begin{remark}
The case $k=2$ is degenerate in \eqref{eq:continuation}
($b\to-\infty$); there $Z_2=\zeta$ directly.
\end{remark}

\newpage

\section{Poles and residues}\label{sec:poles}
\begin{theorem}[Poles and residues]\label{thm:poles}
The only possible poles of $Z_k$ are simple and are located at the half-integers
\begin{equation}\label{eq:pole-locations}
  s=\tfrac12-m,\qquad m=0,1,2,\dots,
\end{equation}
and the residue at $s=\tfrac12-m$ is
\begin{equation}\label{eq:residue}
  \operatorname{Res}_{s=\frac12-m}Z_k
  \;=\;\frac{(-1)^{m}}{2}\binom{2m}{m}\,a^{\,m-\frac12}
  \Bigl(\frac{b^{2}}{16}\Bigr)^{m}.
\end{equation}
\end{theorem}

\begin{proof}
In \eqref{eq:continuation} the only singular part of $\zeta(2s+2j,x)$ is its
pole at $2s+2j=1$, with residue $\tfrac12$ in the variable $s$.
So the pole at $s=\tfrac12-m$ comes from the single term $j=m$:
\[
  \operatorname{Res}=\tfrac12\,a^{-s}\frac{(s)_m}{m!}\Bigl(\frac b2\Bigr)^{2m}
  \Big|_{s=\frac12-m}.
\]
At $s=\tfrac12-m$ we have
$(s)_m=(\tfrac12-m)(\tfrac32-m)\cdots(-\tfrac12)
=(-1)^m\frac{(2m)!}{4^m\,m!}$,
and inserting this gives \eqref{eq:residue}.
\end{proof}

\paragraph{Remarks.}
\begin{itemize}
\item For $m\ge1$ the residue vanishes exactly when $b=0$, i.e.\ $k=4$.
  So $\zeta(2s)$ has the single pole $s=\tfrac12$, and every other
  $Z_k$ has the whole sequence of poles.
\item There are no poles at the negative integers or at $s=0$:
  $Z_k$ is finite there (Section~\ref{sec:trivial}).
\end{itemize}

\begin{example}[Residues]
\begin{center}
\begin{tabular}{@{}lccc@{}}
\toprule
$k$ & residue at $s=\tfrac12$ & residue at $s=-\tfrac12$ & residue at $s=-\tfrac32$\\
\midrule
$3$ & $\dfrac{1}{\sqrt2}$ & $-\dfrac{1}{16\sqrt2}$ & $\dfrac{6}{1024\sqrt2}$\\[1ex]
$4$ & $\dfrac12$ & $0$ & $0$\\[1ex]
$5$ & $\dfrac12\Bigl(\dfrac32\Bigr)^{-1/2}\approx0.408$ & $\approx-0.0085$ & --\\
\bottomrule
\end{tabular}
\end{center}
For $k=3$ the general formula reads
$(-1)^m\binom{2m}{m}\big/\bigl(\sqrt2\cdot32^{m}\bigr)$.
\end{example}

\newpage

\section{Trivial zeros and finite values}\label{sec:trivial}
At non-positive integer arguments the continuation simplifies,
allowing the values of $Z_k$ to be expressed exactly.
At $s=-m$ we have $(s)_j=0$ for $j>m$, so the series \eqref{eq:continuation}
\emph{terminates}:
\begin{equation}\label{eq:terminating}
  Z_k(-m)\;=\;a^{m}\sum_{j=0}^{m}(-1)^{j}\binom{m}{j}
  \Bigl(\frac b2\Bigr)^{2j}\,\zeta(2j-2m,\;x),
\end{equation}
and for non-positive integers $-n$ the Hurwitz value is
\begin{equation}\label{eq:hurwitz-neg}
  \zeta(-n,x)\;=\;-\frac{B_{n+1}(x)}{n+1},
\end{equation}
where $B_r(x)$ is the Bernoulli polynomial.

\begin{proposition}[Value at $s=0$]\label{prop:zero}For every $k>2$,
\begin{equation}\label{eq:Z0}
  Z_k(0)\;=\;-\frac{1}{k-2}.
\end{equation}
\end{proposition}

\begin{proof}
At $s=0$ only $j=0$ contributes, since for $j\ge1$ the factor $(s)_j$ vanishes
and $\zeta(2j,x)$ is finite. So $Z_k(0)=\zeta(0,x)=\tfrac12-x=\tfrac{b-1}{2}$,
and $b-1=-\tfrac{2}{k-2}$.
\end{proof}

\begin{proposition}[Value at $s=-1$]\label{prop:minusone} For every $k>2$,
\begin{equation}\label{eq:Zminus1}
  Z_k(-1)\;=\;\frac{(k-3)(k-4)}{6\,(k-2)^{2}} .
\end{equation}
In particular $s=-1$ is a trivial zero only for $k=3$ and $k=4$.
\end{proposition}

\begin{proof}
Take $m=1$ in \eqref{eq:terminating}:
\[
  Z_k(-1)=a\Bigl[\zeta(-2,x)-\frac{b^{2}}{4}\,\zeta(0,x)\Bigr]
  =a\Bigl[-\frac{B_3(x)}{3}-\frac{b^{2}}{4}\Bigl(\frac12-x\Bigr)\Bigr].
\]
With $B_3(x)=x(x-\tfrac12)(x-1)$ and $x=1-\tfrac b2$ one finds
$-\tfrac13B_3(x)=\tfrac{b(1-b)(2-b)}{24}$ and
$-\tfrac{b^2}{4}(\tfrac12-x)=\tfrac{b^{2}(1-b)}{8}$.
Adding, the bracket equals $\tfrac{b(1-b)(1+b)}{12}$, so
$Z_k(-1)=\tfrac{a\,b\,(1-b^{2})}{12}$.
Finally $ab=\tfrac{k-4}{2}$ and $1-b^{2}=\tfrac{4(k-3)}{(k-2)^{2}}$.
\end{proof}

\begin{theorem}[Trivial zeros of $Z_3$]\label{thm:trivial3}
$Z_3(-m)=0$ for every integer $m\ge1$. 
\end{theorem}

\begin{proof}
For $k=3$ we have $x=\tfrac32$, and from
$B_{2l+1}(\tfrac32)=B_{2l+1}(\tfrac12)+(2l+1)4^{-l}$ with
$B_{\mathrm{odd}}(\tfrac12)=0$ we get
$\zeta(-2l,\tfrac32)=-4^{-l}$.
Then \eqref{eq:terminating} with $a=\tfrac12$, $(b/2)^{2}=\tfrac14$ gives
\[
  Z_3(-m)=-\,2^{-m}\,4^{-m}\sum_{j=0}^{m}(-1)^{j}\binom mj
  =-8^{-m}(1-1)^{m}=0 .
\]
\end{proof}

\begin{remark}
Thus $k=3$ and $k=4$ have trivial zeros at the negative integers. For
other values of $k$, the values $Z_k(-m)$ remain finite. The numerical
examples examined here are generally nonzero, but we do not establish a
general nonvanishing result.
\end{remark}

\newpage

\section{The nontrivial zeros: setting and facts}\label{sec:zeros}
\begin{proposition}[Zero-free half-plane]\label{prop:zerofree}
For $k\ge3$, $Z_k$ has no zeros in $\Re s\ge1$.
\end{proposition}

\begin{proof}
Since $P_k(n)\ge T_n=\tfrac{n(n+1)}2>1$ for $n\ge2$, and
$\sum_{n\ge1}1/T_n=2$ by telescoping, for $\sigma\ge1$
\[
  |Z_k(s)-1|\;\le\;\sum_{n\ge2}P_k(n)^{-\sigma}\;\le\;\sum_{n\ge2}\frac1{T_n}=1,
\]
with strict inequality for $\sigma>1$ or $k>3$.
For $k=3$ and $\sigma=1$, if $Z_3(1+it)=0$, equality in the triangle
inequality would force the terms with $P_3(2)=3$ and $P_3(3)=6$ to have the
same argument as $-1$. Hence $3^{-it}=6^{-it}=-1$. This would make
$\ln3/\ln6$ a ratio of two odd integers, contradicting its irrationality.
The case $t=0$ is excluded by $Z_3(1)=2$.
\end{proof}

For $k \ne 2,4$, no functional equation is available to impose a symmetry
line on the zeros. In what follows, we call the zeros in
$-\tfrac12<\Re s<1$ nontrivial. The numerical computations reported below
place their real parts roughly between $-0.2$ and $0.7$.

\begin{example}[First zeros (numerical)]
\begin{center}
\begin{tabular}{@{}lll@{}}
\toprule
$k=3$ & $k=5$ & $k=7$\\
\midrule
$0.4086+8.8239\,i$  & $0.1502+6.1348\,i$ & $0.0332+5.1129\,i$\\
$0.2085+12.7303\,i$ & $0.2887+9.3010\,i$ & $0.2849+7.9148\,i$\\
$0.5293+15.1042\,i$ & $0.4525+13.5626\,i$ & $0.3603+11.6722\,i$\\
$-0.0190+17.7999\,i$ & $-0.0528+15.1236\,i$ & $0.0792+13.7607\,i$\\
\bottomrule
\end{tabular}
\end{center}
Of course, zeros come in complex conjugate pairs.
\end{example}

\begin{definition}[Window mean]\label{def:window}
Fix $k>2$, a height $t>0$ and a window length $h>0$.
Let $\mathcal Z$ be the set of zeros $\rho$ of $Z_k$ (counted with multiplicity)
with
\[
  - \tfrac 12 < \Re\rho < 1,\qquad t<\Im\rho\le t+h,
\]
let $N=|\mathcal Z|$, and define the arithmetic mean of the real parts
\begin{equation}\label{eq:window-mean}
  \bar\sigma_k(t,h)\;=\;\frac1N\sum_{\rho\in\mathcal Z}\Re\rho .
\end{equation}
The centre of the window is $t_c=t+h/2$. The zero-free result in
Proposition~
ef{prop:zerofree} justifies the right boundary at $1$.
The left boundary at $-	\tfrac12$ is a convention; the computations so far
have not produced a zero with real part below approximately $-0.2$.
\end{definition}

The law below needs windows that are long but not too long:
\begin{equation}\label{eq:window-condition}
  \ln t\;\ll\;h\;\ll\;t .
\end{equation}
This is a working condition from error estimates, not a proved requirement.

\newpage
\section{Derivation of the mean real part law}\label{sec:law}
\textbf{Notation.}
\begin{equation}\label{eq:notation-law}
  \ell=\ln\frac{t_c}{\pi},\qquad
  a=\frac{k-2}{2},\qquad
  \Lambda=2\ell+\ln a=\ln\frac{(k-2)\,t_c^{2}}{2\pi^{2}} .
\end{equation}
\subsection*{The zero-free region to the right}
The first term of $Z_k$ is exactly $1$ (Eq.~\eqref{eq:first-term}), and for
$\sigma>1$ the remaining terms are smaller than $1$ in total, so $Z_k\ne0$
(Proposition~\ref{prop:zerofree}).

\subsection*{A mean-value function on vertical lines}
For a fixed real part $\sigma$, define the mean of $\log|Z_k|$ over the
height interval by
\begin{equation}\label{eq:M}
  M(\sigma)\;=\;\frac1h\int_t^{t+h}\ln\bigl|Z_k(\sigma+iu)\bigr|\,du .
\end{equation}
for fixed $\sigma>1$, the absolutely convergent
Dirichlet-series expansion of $\log Z_k(s)=\log(1+w(s))$ has only positive
frequencies. Each individual frequency has a nonzero finite-window average,
but its average is $O(1/h)$. Absolute convergence permits termwise
averaging, and therefore
\begin{equation}\label{eq:M-right}
  M(\sigma)=O_{k,\sigma}(1/h)\qquad(\sigma>1).
\end{equation}
Thus this mean tends to zero as the window length grows, although it need
not vanish for any fixed finite window.

\subsection*{Application of Littlewood's lemma}
\textbf{[Standard theorem]}\;
Littlewood's lemma (the strip version of Jensen's formula; see
Titchmarsh) relates the change in this mean to the zeros inside a rectangle.
Each zero contributes according to its horizontal displacement from the
left boundary.
Let $\nu=N/h$ be the number of zeros per unit height and $\bar\sigma$ their
mean real part.
Take $-0.5\leq\sigma_0\leq-0.4$ and a right boundary $\sigma_1>1$.
Under Hypothesis~\ref{hyp:H}, there are no zeros to the left of this interval;
Proposition~\ref{prop:zerofree} excludes zeros to the right of $1$. Thus the
zeros in the rectangle are precisely those in the window definition, provided
there are no zeros on its horizontal edges (otherwise move the edges
slightly). Littlewood\textquotesingle s lemma, applied to this rectangle,
gives
\begin{equation}\label{eq:jensen}
  M(\sigma_0)\;=\;2\pi\nu\,\bigl(\bar\sigma-\sigma_0\bigr)\;+\;O\!\left(\frac{\ln t}{h}\right),
\end{equation}
provided the horizontal-edge contribution has the stated bound. This
horizontal-edge estimate is an additional assumption, uniform for
$-0.5\leq\sigma_0\leq-0.4$; it is not a consequence of Littlewood\textquotesingle s lemma alone.

\emph{Interpretation.} To the left of the zeros, the formula describes a
linear function of $\sigma_0$. Its
\begin{itemize}
\item \textbf{slope} is $-2\pi\nu$ (it tells us how many zeros there are), and
\item \textbf{value at $\sigma_0=0$} is $2\pi\nu\,\bar\sigma$ (it tells us how far to the right
  the zeros are, weighted by their number).
\end{itemize}

\subsection*{The mean value to the left of the zeros}
Far to the left, $Z_k=a^{-s}H_b(s)$ is dominated by a single ``mirror'' term,
exactly as for the Riemann zeta function, whose size is
$(t/\pi)^{1/2-2\sigma}$, while the prefactor $a^{-s}$ has size $a^{-\sigma}$.
Taking logarithms and averaging motivates the following hypothesis; the
existence and dominance of such a dual series are not established here.

\begin{hypothesis}[Mean value to the left of the zeros]\label{hyp:H}
\textbf{[Assumed]}\;
There are no zeros with $\Re\rho<-0.4$, and for $-0.5 \le \sigma_0 \le -0.4$,
\begin{equation}\label{eq:H}
  M(\sigma_0)\;=\;\Bigl(\tfrac12-2\sigma_0\Bigr)\ell\;-\;\sigma_0\ln a\;+\;o(1).
\end{equation}
\end{hypothesis}

\emph{Where it comes from.}
More specifically, a Riemann--Siegel-type formula of the form
$H_b(s)=\chi(s)\widetilde H(s)$, with
$|\chi(s)|=(t/\pi)^{1/2-2\sigma}$, would lead to the stated expression if
the leading term of the dual series $\widetilde H$ dominated for
$\sigma\le-0.36$ and its logarithmic modulus had mean zero. The additional
term $-\sigma_0\ln a$ is contributed by the factor $a^{-s}$.

\emph{Status.}
For $k=4$ this is classical (it is the functional equation of $\zeta(2s)$).
For other values of $k$, the hypothesis remains unproved. The table
below compares its prediction with the reported numerical mean.

\begin{table}[h]
\centering
\caption{Direct test of Hypothesis~\ref{hyp:H}: average of $\ln|Z_k(-0.45+iu)|$
for $1000<u<1300$ against the right-hand side of \eqref{eq:H}.
\textbf{[Numerical]}}\label{tab:H}
\begin{tabular}{@{}lccccc@{}}
\toprule
$k$ & $3$ & $5$ & $11$ & $101$ & $1001$\\
\midrule
measured $M(-0.45)$  & $7.948$ & $8.443$ & $8.937$ & $10.016$ & $11.056$\\
predicted            & $7.952$ & $8.446$ & $8.941$ & $10.020$ & $11.060$\\
\bottomrule
\end{tabular}
\end{table}

\subsection*{Comparison of the two mean-value formulas}
Rewrite \eqref{eq:H} as a straight line in $\sigma_0$:
\begin{equation}\label{eq:line}
  M(\sigma_0)\;=\;\frac{\ell}{2}\;-\;\sigma_0\,(2\ell+\ln a)
  \;=\;\frac{\ell}{2}-\sigma_0\Lambda .
\end{equation}
Compare this with \eqref{eq:jensen}. If the $o(1)$ in Hypothesis~\ref{hyp:H}
and the horizontal-edge error are uniform on the interval
$-0.5\leq\sigma_0\leq-0.4$, evaluating at two distinct fixed values of
$\sigma_0$ gives agreement of the slopes and intercepts up to $o(1)$:
\begin{align}
  \text{slope:}\quad & 2\pi\nu=\Lambda+o(1), \label{eq:slope}\\
  \text{intercept:}\quad & 2\pi\nu\,\bar\sigma=\frac{\ell}{2}+o(1). \label{eq:intercept}
\end{align}

\newpage

\section{The mean real part law}

\begin{hypothesis}[The law]\label{thm:law}
\textbf{[Conditional on Hypothesis~\ref{hyp:H}, uniformity of its error,
and the horizontal-edge estimate in Step~3.]}; In a window with centre $t_c$
satisfying the window condition
\eqref{eq:window-condition},
the number of zeros per unit height and their mean real part satisfy
\begin{equation}\label{eq:law}
  \boxed{\;
  \nu_k=\frac{\Lambda}{2\pi}+o(1),
  \qquad
  \bar\sigma_k\;=\;\frac{\ell}{2\Lambda}
  \;+\;o\!\left(\frac1{\Lambda}\right)
  \;=\;\frac14-\frac{\ln a}{4\Lambda}
  \;+\;o\!\left(\frac1{\Lambda}\right)\;}
\end{equation}
\end{hypothesis}

\begin{proof}
Dividing the intercept \eqref{eq:intercept} by the slope \eqref{eq:slope}
gives $\bar\sigma=\ell/(2\Lambda)+o(1/\Lambda)$, since
$\Lambda=2\ell+\ln a$ grows like $2\ell$ for fixed $k$ as $t_c\to\infty$.
\end{proof}

\begin{corollary}[Consequences of the law]\label{cor:limits}
\begin{itemize}
\item $k=4$ ($a=1$): the functional equation and the symmetry of the
  zeros of $\zeta(2s)$ about $\Re s=\tfrac14$ give mean real part
  $\tfrac14$ in every height window, by the symmetry pairing zeros with
  real parts $\beta$ and $1-\beta$ at the same height; their density is
  asymptotic to $\tfrac1\pi\ln\tfrac t\pi$.
\item In the regime $\Lambda>0$, $k<4$ ($a<1$) predicts a mean
  \emph{above} $\tfrac14$, while $k>4$ ($a>1$) predicts a mean
  \emph{below} $\tfrac 14$.
\item Fixed $k$, height $t\to\infty$: $\bar\sigma_k\to\tfrac14$.
\item If the law were uniform in $k$ at fixed height, it would predict
  $\bar\sigma_k\approx\ell/(2\ln a)\to0$, like $1/\ln k$.
  The fixed-$k$ asymptotic proved conditionally above does not by itself
  justify this interchange of limits. 
\item As $k\to2^+$ at fixed height, $\ln a\to-\infty$ and the displayed
  expression eventually leaves the regime $\Lambda>0$ required for a positive
  zero density. The law therefore gives no conclusion in this limit; in
  particular it does not imply convergence to $\tfrac12$. 
\end{itemize}
\end{corollary}

\newpage
\section{The law in use: examples and numerical checks}\label{sec:examples}
All data below are \textbf{[Numerical]}: zeros were located with a
Euler--Maclaurin evaluator of $Z_k$, and in the windows used the number of
zeros found matched the argument-principle count, except where a shortfall is
stated.

\begin{example}[$k=5$ at height $10\,000$]
Here $a=1.5$, $\ell=\ln(10000/\pi)=8.07$, $\ln a=0.405$, so $\Lambda=16.54$.
The law predicts a mean $8.07/(2\cdot16.54)\approx0.244$ and
$\Lambda/2\pi\approx2.6$ zeros per unit height.
Measured: mean $0.2475\pm0.0037$ in a window of length $400$.
\end{example}

\begin{table}[h]
\centering
\caption{Number of zeros in windows of length 30 at $t\approx1015$:
counted via the argument principle against $\Lambda h/(2\pi)$.}
\begin{tabular}{@{}lcccccccc@{}}
\toprule
$k$ & 3 & 5 & 7 & 11 & 21 & 51 & 101 & 1001\\
\midrule
counted & 52 & 58 & 60 & 62 & 66 & 70 & 74 & 85\\
formula & 51.9 & 57.1 & 59.6 & 62.4 & 65.9 & 70.4 & 73.8 & 84.8\\
\bottomrule
\end{tabular}
\end{table}

\begin{table}[h]
\centering
\caption{Mean real part of the zeros: observed against the law
\eqref{eq:law}. ($t_c$ is the window centre.)}
\begin{tabular}{@{}rrcc@{}}
\toprule
$t_c$ & $k$ & observed & law\\
\midrule
120   & 3    & 0.275 & 0.276\\
120   & 4    & 0.250 & 0.250\\
120   & 5    & 0.234 & 0.237\\
120   & 11   & 0.211 & 0.207\\
120   & 101  & 0.164 & 0.163\\
120   & 1001 & 0.139 & 0.135\\
120   & $10^{5}+1$ & 0.103 & 0.101\\
120   & $10^{7}+1$ & 0.082 & 0.080\\
1015  & 3    & 0.262 & 0.266\\
1015  & 5    & 0.250 & 0.242\\
1015  & 101  & 0.186 & 0.187\\
1015  & 1001 & 0.164 & 0.163\\
10015 & 11   & 0.236 & 0.229\\
10015 & 101  & 0.203 & 0.201\\
\bottomrule
\end{tabular}
\end{table}

Over all $25$ comparisons ($k$ from $2.05$ to $10^7+1$, $t_c$ from $120$ to
$10^4$) the root-mean-square error is $0.008$; for $k\ge2.2$ no error exceeds
$0.008$.
The law overshoots at $k=2.1$ ($+0.011$) and $k=2.05$ ($+0.033$), where
$a\to0$ and it is expected to fail.

\begin{table}[h]
\centering
\caption{$k=5$ at large height: arithmetic mean $\pm$ standard error against
the law (found/expected = zeros located / $\Lambda h/2\pi$).}
\begin{tabular}{@{}rrcccc@{}}
\toprule
$t_c$ & $h$ & found/expected & mean & law & difference/s.e.\\
\midrule
450    & 300 & 491/493  & $0.2408\pm0.0067$ & 0.2402 & $+0.1$\\
2100   & 200 & 425/427  & $0.2438\pm0.0065$ & 0.2424 & $+0.2$\\
5075   & 150 & 359/362  & $0.2461\pm0.0067$ & 0.2433 & $+0.4$\\
10030  & 60  & 156/157  & $0.2445\pm0.0098$ & 0.2439 & $+0.1$\\
10200  & 400 & 1044/1055 & $0.2475\pm0.0037$ & 0.2439 & $+1.0$\\
30075  & 150 & 443/447  & $0.2474\pm0.0056$ & 0.2446 & $+0.5$\\
100030 & 60  & 196/202  & $0.2551\pm0.0071$ & 0.2452 & $+1.4$\\
\bottomrule
\end{tabular}
\end{table}

The $k=5$ means rise toward $\tfrac14$ \emph{from below}, as the law predicts,
while for $k=3$ the means fall toward $\tfrac14$ \emph{from above}:

\begin{table}[h]
\centering
\caption{$k=3$: observed mean against the law.}
\begin{tabular}{@{}rccc@{}}
\toprule
$t_c$ & observed & s.e. & law\\
\midrule
500    & 0.2716 & 0.0074 & 0.2683\\
1150   & 0.2667 & 0.0073 & 0.2656\\
3075   & 0.2627 & 0.0091 & 0.2632\\
10050  & 0.2620 & 0.0094 & 0.2612\\
30030  & 0.2650 & 0.0105 & 0.2598\\
100015 & 0.2650 & 0.0123 & 0.2586\\
\bottomrule
\end{tabular}
\end{table}

All seven $k=5$ differences are positive (average about $+0.003$).
This is consistent with the small shortfalls in the counts: missed zeros tend
to have small real part, which biases the mean upward.
In the one window with a verified count ($156$ of $157$) there is no visible
bias.

\newpage

\section{The one-fourth concentration hypothesis}\label{sec:hypothesis}
By our hypothesized law \eqref{eq:law}, for every fixed $k \ge 3$, and $t_c = t + h/2$ satisfying the window condition \eqref{eq:window-condition} we have
\begin{hypothesis}[Mean to Quarter]
  \label{eq:mean-to-quarter}
  \begin{equation}
  \bar\sigma_k(t,h) \;=\;\frac14-\frac{\ln a}{4\Lambda}
  \;+\;o\!\left(\frac1{\Lambda}\right)\;\longrightarrow\;\frac14
  \qquad(t_c\to\infty) 
  \end{equation}
\end{hypothesis}
This conjecture concerns the mean real part only. Even if the mean tends
to $1/4$, it does not follow that the zeros themselves concentrate near that
line: a small number of zeros far from $1/4$ could have a substantial effect
on the arithmetic mean. The next hypothesis is therefore strictly stronger.
It asserts that the proportion of zeros lying outside any fixed
neighbourhood of $\Re s= \tfrac14$ tends to zero as the height increases.
The available computations cover only a modest range of heights and cannot
determine whether such concentration occurs asymptotically. 

\begin{hypothesis}[One-fourth Concentration]\label{hyp:quarter}
Fix $k \ge 3$ and let $X$ be the real part of a zero chosen uniformly among the
zeros in a window of centre $t_c$ and length $h$ obeying the window condition
\eqref{eq:window-condition}.
Then for every $\delta>0$
\begin{equation}\label{eq:concentration}
  \mathbb P\bigl(\,|X-\tfrac14|>\delta\,\bigr)\;\longrightarrow\;0
  \qquad(t_c\to\infty).
\end{equation}
\end{hypothesis}

\paragraph{Numerical evidence for the shrinking.}
The standard deviation of the real parts falls with height:
\begin{center}
\begin{tabular}{@{}lccccccc@{}}
\toprule
$t_c$ (for $k=3$) & 120 & 500 & 1150 & 3075 & 10050 & 30030 & 100015\\
s.d. & 0.197 & 0.179 & 0.167 & 0.161 & 0.147 & 0.135 & 0.118\\
\midrule
$t_c$ (for $k=5$) & 450 & 2100 & 5075 & 10200 & 30075 & 100030 & \\
s.d. & 0.149 & 0.134 & 0.128 & 0.121 & 0.117 & 0.100 & \\
\bottomrule
\end{tabular}
\end{center}
The standard deviations decrease across the ranges shown. This trend is
consistent with a narrowing distribution, but finite numerical data do not
establish concentration in the limit.

\paragraph{What the evidence does not show.}
\begin{itemize}
\item The product $\mathrm{s.d.}\times\Lambda$ is \emph{not} constant: for $k=5$
  it grows slowly from about $1.5$ to $2.1$, and similarly for $k=3$.
\item At $t_c\ge10^4$ the standard errors of the means are $0.004$ to $0.012$,
  while the offset of the law from $\tfrac14$ is only $0.005$ to $0.015$.
  The data therefore separate the law from a plateau far from $\tfrac14$,
  but they cannot separate the law from the limit $\tfrac14$ itself.
\item The distribution is skewed, with a longer right tail; the median is
  slightly smaller than the arithmetic mean. The mean should therefore be
  interpreted with some care.
\end{itemize}

\paragraph{What would turn the hypotheses into theorems.}
\begin{enumerate}
\item A rigorous approximate functional equation for $H_b(s)$, with the
  dual leading coefficient of modulus $1$ and error $o(1)$ (Poisson summation
  and a van der Corput $B$-process for the quadratic phase).
\item A bound on the horizontal-edge term in Littlewood's lemma for $Z_k$.
\item A tent estimate, i.e.\ control of the mean of
  $\ln|Z_k|$ for $\sigma$ between the zeros and the right half-plane.
\item A proof that $Z_k$ actually has the predicted number of zeros
  (infinitely many) for $k\ne2,4$.
\end{enumerate}

\newpage

\section{The limit $k\to\infty$: Do the zeros wander onto the imaginary axis?}\label{sec:limit}
We now consider the behaviour of the zeros as $k$ increases. The
asymptotic law suggests that the mean real part of the nontrivial zeros may
move towards the imaginary axis, but the limiting behaviour requires care.
In particular, a separate zero approaches $s=\tfrac12$. We write
$L=k-2$ throughout this section.

\paragraph{Why $Z_\infty$ does not exist as a function.}
For $\Re s>0$ the limit of $Z_k(s)$ away from poles as $k\to\infty$ is the constant $1$.
The decomposition proved below is restricted to $\Re s>0$ and gives no
conclusion about the size of $Z_k(s)$ when $\Re s<0$. Thus the limit on the
whole plane is not established, and there is no limit function on that domain
to which one can assign zeros.
The meaningful question is where the zeros of $Z_k$ accumulate as $k\to\infty$.

\begin{proposition}[Exact decomposition for $\Re s>0$]\label{prop:decomposition}
For $L\ge2$ and $\Re s>0$,
\begin{equation}\label{eq:decomposition}
  Z_k(s)\;=\;1+L^{-s}\,Z_3(s)+\Bigl(\frac L2\Bigr)^{-s}\Delta_k(s),
\end{equation}
where
\begin{equation}\label{eq:Delta}
  \Delta_k(s)=\sum_{m\ge1}\bigl(m(m+1)\bigr)^{-s}
  \Bigl[\Bigl(1+\frac{2}{Lm}\Bigr)^{-s}-1\Bigr]
\end{equation}
converges absolutely and
\begin{equation}\label{eq:Delta-bound}
  |\Delta_k(s)|\;\le\;\frac{2|s|}{L}\,\zeta(1+2\sigma).
\end{equation}
\end{proposition}

\begin{proof}
For $\sigma>\tfrac12$ the terms with $n\ge2$ of $Z_k$ are
$P_k(n)^{-s}$ with $P_k(n)=\tfrac L2\,n\bigl(n-1+\tfrac2L\bigr)$.
Put $m=n-1$ to get
$\sum_{n\ge2}P_k(n)^{-s}=(L/2)^{-s}\sum_{m\ge1}\bigl(m(m+1)\bigr)^{-s}
\bigl(1+\tfrac{2}{Lm}\bigr)^{-s}$.
Split $\bigl(1+\tfrac2{Lm}\bigr)^{-s}=1+\bigl[\bigl(1+\tfrac2{Lm}\bigr)^{-s}-1\bigr]$.
The first part gives $(L/2)^{-s}\sum_m(m(m+1))^{-s}=(L/2)^{-s}2^{-s}Z_3(s)
=L^{-s}Z_3(s)$, since $T_m=m(m+1)/2$.
The second part is $(L/2)^{-s}\Delta_k(s)$.
For $0\le y\le1$ and $\sigma\ge0$,
$\bigl|(1+y)^{-s}-1\bigr|=\bigl|\int_0^y s(1+u)^{-s-1}du\bigr|\le|s|\,y$,
and $\bigl(m(m+1)\bigr)^{-\sigma}\le m^{-2\sigma}$, so each term of
\eqref{eq:Delta} is at most $\frac{2|s|}{L}m^{-1-2\sigma}$ in modulus.
This proves \eqref{eq:Delta-bound} and shows $\Delta_k$ is analytic for
$\sigma>0$.
Both sides of \eqref{eq:decomposition} are meromorphic on $\sigma>0$ and agree
for $\sigma>\tfrac12$, so they agree on $\sigma>0$ by the identity theorem.
\end{proof}

\begin{theorem}[Zeros leave every compact set in $\Re s>0$ except near $\tfrac12$]
\label{thm:expelled}
Let $K\subset\{\Re s>0\}$ be compact with $\tfrac12\notin K$.
Then $\sup_{s\in K}|Z_k(s)-1|\to0$ as $k\to\infty$.
In particular there is $k_0(K)$ such that $Z_k$ has no zeros in $K$ for
$k\ge k_0$.
\end{theorem}

\begin{proof}
On $K$ we have $\sigma\ge\sigma_K>0$, $|s|$ is bounded, and $Z_3$ is bounded
(its only pole in $\Re s>0$ is at $s=\tfrac12$).
By \eqref{eq:decomposition} and \eqref{eq:Delta-bound},
\[
  |Z_k(s)-1|\;\le\;L^{-\sigma_K}\sup_K|Z_3|
  +\Bigl(\frac L2\Bigr)^{-\sigma_K}\frac{2\sup_K|s|}{L}\,\zeta(1+2\sigma_K)
  \;\longrightarrow\;0 .
\]
\end{proof}

\begin{theorem}[One real zero near the pole at $\tfrac12$]\label{thm:realzero}
For all sufficiently large $k$, the disc $|s-\tfrac12|\le2L^{-1/2}$ contains
exactly one zero of $Z_k$; it is real, and
\begin{equation}\label{eq:realzero}
  s_*\;=\;\frac12-\frac{1}{\sqrt{2L}}\;\bigl(1+o(1)\bigr).
\end{equation}
\end{theorem}

\begin{proof}
Write $s=\tfrac12+L^{-1/2}z$ with $|z|\le2$, and
$Z_3(s)=\dfrac{1}{\sqrt2\,(s-\frac12)}+g(s)$, where $g$ is analytic and bounded
near $\tfrac12$ (the residue $\tfrac1{\sqrt2}$ is from Section~\ref{sec:poles}).
By \eqref{eq:decomposition},
\[
  z\,Z_k(s)\;=\;z+\frac{1}{\sqrt2}\,e^{-zL^{-1/2}\ln L}
  +z\,L^{-s}g(s)+z\,(L/2)^{-s}\Delta_k(s).
\]
As $L\to\infty$ the second term tends to $\tfrac1{\sqrt2}$ uniformly on
$|z|\le2$, because $L^{-1/2}\ln L\to0$, and the last two terms are
$O(L^{-1/2})$.
So $z\,Z_k\to z+\tfrac1{\sqrt2}$ uniformly on $|z|\le2$.
This limit has the single zero $z=-\tfrac1{\sqrt2}$, and the left side does not
vanish at $z=0$ (its value there tends to $\tfrac1{\sqrt2}$), so the zeros of
$zZ_k$ in the disc are exactly those of $Z_k$.
By Rouch\'e's theorem there is exactly one, near $z=-\tfrac1{\sqrt2}$.
Since $Z_k(\bar s)=\overline{Z_k(s)}$ and the disc is symmetric, the zero is
real.
\end{proof}

\paragraph{Consequence.}
The point $s=\tfrac12$ is a limit point of zeros of $Z_k$, and it does
\emph{not} lie on the imaginary axis.
Thus it is not correct to say that all zeros approach the imaginary axis.
What the theorem proves is that zeros cannot accumulate in the open right
half-plane away from $s= \tfrac12$.

\begin{table}[h]
\centering
\caption{Numerical check of Theorem~\ref{thm:realzero}: the real zero of $Z_k$
in $(0,\tfrac12)$ against the prediction $\tfrac12-1/\sqrt{2L}$.
\textbf{[Numerical]}}
\begin{tabular}{@{}rccc@{}}
\toprule
$k$ & $11$ & $101$ & $1001$ and $10^5+1$\\
\midrule
real zero found & none & $0.3918$ & $0.4737$,\ $0.49771$\\
prediction      & $0.2643$ & $0.4289$ & $0.4776$,\ $0.49776$\\
\bottomrule
\end{tabular}
\end{table}

At $k=11$ ($L=9$), the asymptotic regime is not yet effective. The
discrepancy is smaller at $k=101$ and decreases as $k$ increases.

\paragraph{A second real zero near $0$.}
For $k\ge101$ the numerics also show a tiny positive real zero
($0.0049$, $0.00023$ and about $0$ for $k=101$, $1001$, $10^5+1$).

\paragraph{The remaining zeros: conjecture and evidence.}
The results in this section leave the remaining zeros largely unresolved.
The following conclusions are conditional or numerical, rather than
theorems:
\begin{itemize}
\item Conditional on Hypothesis~\ref{hyp:H}, the zero density is
  $\Lambda/(2\pi)\to\infty$ and the mean real part is
  $\ell/(2\ln a)\to0$ (Corollary~\ref{cor:limits}), so the zeros become
  infinitely dense and their mean position tends to the imaginary axis.
  This is a statement about the average, not about every zero.
\item The numerics show that the zeros are \emph{not literally on} the axis:
  for $k=101$, $1001$, $10^5+1$ there are zeros at low height with
  \emph{negative} real part, for example
  $-0.202+2.034\,i$ ($k=101$),
  $-0.149+1.250\,i$ and $-0.114+2.290\,i$ ($k=1001$),
  $-0.061+0.674\,i$ and $-0.083+1.302\,i$ ($k=10^5+1$).
  They move toward the axis only slowly.
  Nothing is proved in $\Re s<0$.
\end{itemize}

\newpage

\section*{Declaration of the Use of Artificial Intelligence (AI)}

During the preparation of this research paper, Artificial 
Intelligence (AI) tools and referees, including Claude Sonnet 5.5 as a main assistant and GPT 6.0 as a referee, were 
used as supplementary aids to support the research and 
writing process. Specifically, AI assistance was utilized 
for improving language clarity, grammar, spelling, sentence 
structure, academic phrasing, and the overall readability 
of the manuscript. AI tools also assisted in organizing 
ideas, refining expressions, summarizing complex concepts, 
and providing suggestions for improving the structure 
and coherence of the paper.

The use of AI was limited to supporting and enhancing 
the writing process and did not replace the researcher's 
intellectual contribution, critical thinking, independent 
analysis, interpretation of findings, or academic judgment. 
The research objectives, methodology, data collection, 
analysis, interpretation, and conclusions remain the 
responsibility of the author(s).

All AI-generated suggestions and outputs were critically 
reviewed, evaluated, and revised by the author(s) to ensure 
accuracy, relevance, originality, and consistency with the 
research objectives. The author(s) take full responsibility 
for the content, accuracy, integrity, and final version of 
this manuscript, including the verification of all references, 
sources, facts, and research findings.

The author(s) acknowledge the importance of transparency, 
academic integrity, and the responsible use of emerging 
AI technologies in scholarly research. AI tools were not 
considered authors or contributors in the academic authorship 
sense, and their use does not diminish the independent 
scholarly contribution of the author(s).

\newpage

\section*{Summary: what is proved and what is assumed}
\begin{center}
\begin{tabular}{@{}p{8.2cm}p{5.6cm}@{}}
\toprule
Statement & Status\\
\midrule
Continuation, poles, residues, $Z_k(0)$, $Z_k(-1)$ & Proved\\
Trivial zeros of $Z_3$ at $-1,-2,\dots$ & Proved\\
No zeros in $\Re s\ge1$ & Proved\\
Right-hand mean $M(\sigma)=O(1/h)$ for $\sigma>1$ & Proved\\
Jensen's identity on a strip & Standard; edge error assumed small\\
Mean value to the left of the zeros (Hypothesis~\ref{hyp:H}) & Assumed; classical for $k=4$; numerically checked to $0.004$\\
Density $\Lambda/2\pi$ and mean $\ell/2\Lambda$ & Conditional on the above; numerical error about $1\%$ and $0.008$\\
Concentration at $\tfrac14$ & Hypothesis; needs a tent estimate\\
For $k$ large: no zeros in compacts of $\Re s>0\setminus\{\tfrac12\}$; one real zero near $\tfrac12-1/\sqrt{2(k-2)}$ & Proved\\
\bottomrule
\end{tabular}
\end{center}

\end{document}